\documentclass[12pt,twoside, epsf]{article}

\usepackage{color,graphicx,times}

\makeatother

\let \ttorg \tt \def \tt{\ttorg \obeyspaces}

\begin{document}

\title{\bf  Book Review of  ``Virtual Knot Theory - The State of the Art" 
by Vassily Manturov and Denis Ilyutko, World Scientific (2012), ISBN: 978-981-4401-12-8}

\author{Louis H. Kauffman\\ Department of Mathematics, Statistics \\ and Computer Science (m/c
249)    \\ 851 South Morgan Street   \\ University of Illinois at Chicago\\
Chicago, Illinois 60607-7045\\ $<$kauffman@uic.edu$>$\\}

\maketitle

\thispagestyle{empty}

This book, is the first full-length book on the theory of virtual knots and links. I coined the term 
virtual knot theory and wrote a first paper \cite{VKT} on this subject in 1999. From my point of view, the idea 
was to simultaneously have a diagrammatic theory that could handle knots in thickened surfaces and 
would generalize knot theory to arbitrary oriented, not necessarily planar, Gauss codes. These 
non-planar Gauss codes act in knot theory as abstract graphs behave in graph theory. They can be projected into the plane, and will acquire virtual crossings in the process. They can be embedded in surfaces of
higher genus, and one can search for the surface of least genus that will support them. The abstract Gauss codes (often depicted as so-called Gauss diagrams) have Reidemeister moves defined for them
irrespective of a choice of embedding of the code in a surface. Once embedded, they yield a knot in
a thickened version of that surface. Classical knot theory embeds in virtual knot theory.

Many invariants of classical knots generalize directly to virtual knots,
the fundamental group and the quandle among them. The first quantum invariant that has an easy generalization is the Jones polynomial. The bracket state sum model of the Jones polynomial generalizes directly to virtual knot theory and has a number of startling properties. One of these 
properties is that a non-trivial classical knot can be transformed to a virtual knot by changing only the orientation of some of its crossings (we call this process the {\it virtualization of the crossing}).
If a subset of crossings is chosen so that switching those crossings unknots the classical knot, then the
virtual knot obtained by reversing their orientations will be a non-trivial virtual knot with Jones polynomial equal to one. Thus there are infinitely many non-trivial virtual knots with unit Jones polynomial. We would like to know if any of these are actually classical knots. We suspect that this is
not the case and the research on this problem has led to the invention/discovery of new invariants by 
myself and by many other workers in virtual knot theory, including particularly the authors of this book and Vassily Manturov most particularly. This problem about the Jones polynomial is a driving force behind much of the research on virtual knot theory. Of course, in back of that problem is the question 
whether the Jones polynomial detects the unknot for classical knots.

All quantum invariants of classical knots generalize to a special category of virtual knots and links
that I call the {\it rotational virtuals} and are called in this book the {\it rigid virtuals}. Using the present 
terminology, a rigid virtual knot does not have a flat version of the first Reidemeister move available for its diagrams. In \cite{VKT} I began the study of this version of virtual knot theory and it is still in need of 
development. The reader will find a short exposition of rigid virtual knot theory in the present book in
Chapter 4, in the paper \cite{VKT} and in the paper \cite{VBCAT} where we relate this theory both to a categorical 
approach to the virtual braid group and to quantum invariants derived from Hopf algebras (using the terminology ``rotational"). In any case, 
this underlines the meaning of  ``State of the Art" in the title of the present book. There is much that is 
in the course of development in virtual knot theory, and this book is a snapshot of that development from 
the point of view of the authors of the book.

Turning to the book itself, the reader will find a remarkable and encyclopedic treatment of the basics, with the first chapter devoted to the diagrammatic and combinatorial definitions and to a discovery of
mine of which I am very fond -- the self-linking number of a virtual knot \cite{SL}. This self-linking number is,
in the Gauss-diagram language the sum of the crossing signs for those chords in the diagram that are
{\it odd}. A chord in a chord diagram is odd if it is intersected by an odd number of other chords in the diagram. The appearance of odd chords is a feature that can only happen in a virtual knot or link. The self-linking number is surprisingly powerful and certainly the simplest combinatorial invariant that one can find for virtual knots. It is also the tip of the iceberg of {\it parity}, a subject that Manturov has expanded for virtual knot theory and found ways to exploit that are highly significant for virtual theory and, we 
expect for classical theory as well in an indirect way. Later in the book, in Chapter 8, this full development of Maturov Parity Theory for Virtual Knots is given in its entirety. Here again, we go from the 
basics of the first chapter to the frontiers of research. The surface structure of parity that one sees in the 
Gauss diagrams for virtual knots and links is certainly just an indication of deeper matters of parity that 
are in the works below that surface. The future of virtual knots and their relationship with classical knots
turns on this deep parity.

Continuing to Chapter 2, we find a discussion of Kuperberg's Theorem \cite{Kup} and the surface genus of 
virtual knots. Kuperberg proved the remarkable theorem that if a virtual knot or link is represented 
in its minimal genus surface, then this embedding type is unique. This means that the virtual knots 
have a definite topological intepretation  in their minimal genus surfaces. This result of Kuperberg
has been very fruitful for getting deeper invariants of virtual knots. Heather Dye and myself found 
ways to use the Kuperberg theorem to get stronger version of the Jones polynonmial for virtuals
\cite{MinSurf,Arrow} and Manturov did a similar analysis that will be found in this book. This chapter proves that
virtual knots are algorithmically recognizable by generalizing the technique of Haken and Hemion
for classical knots and it proves that connected sums of non-trivial virtual knots are non-trivial.
This is particularly significant in the light of the phenomenon that connected sums of trivial virtual diagrams can yield non-trivial virtual knots!

In Chapter 3 the authors consider quandles and their generalizations for virtual knots.
One of the most fruitful generalizations of invariants that has emerged from virtual knot theory is the use of the biquandle \cite{BIQ,VBIQ}, a generalization of the quandle (which in turn generalizes the fundamental group of a knot or link) and its powerful applications to the theory. Here also the reader will find Manturov's remarkable uses of Lie algebraic techniques for invariants. In this same chapter one will find long 
virtual knots and clever invariants for them using Manturov's application of the precendence order that
is available there and flat virtuals where one is looking at immersions of curves in surfaces. This leads to
a notion of mine called the (flat) hierarchy where we have as many virtual crossings as one likes, ordered by some choice of an ordinal, in terms of their ability to move across one another. Maturov's first steps at  some invariants are found here for the hierarchy.

Chapter 4 does the basics of the Jones polynomial via the bracket state sum and introduces the
concept of atoms (surfaces, orientable or not, bearing the virtual knot diagram). Now the pace picks up as we go to Chapter 5 and meet Khovanov homology for virtual knots. Manturov gave a first combinatorial solution to the question of constructing integral Khovanov homology for virtual knots and 
links. This has led to initial collaborations on virtual Khovanov homology \cite{Categor,Kaestner}. His theory is quite interesting and it needs to be compared with the Khovanov-Rozansky 
Categorified Link Homology and it needs to be computed. There is a rich vein of research material for the active reader who gets to this point in the book.

We are not done! Chapter 6 deals with virtual braids and the work of Bardakov and of Kamada
and Kauffman and Lambropoulou \cite{VBL}. It ends with an exposition of Manturov's invariants of virtual braids.
Chapter 7 treats Vassiliev invariants and returns to another source of virtual knot theory -- the work of 
Goussarov, Polyak and Viro \cite{GPV} who used virtual knots in the guise of general Gauss diagrams to construct
a theory of Gauss diagram formulas for virtual knots.  Here we find that theory and a exposition of 
Manturov's proof of a conjecture of Vassiliev about these diagrams. Then comes Chapter 8 on Parity, alluded to before. This Chapter includes work on the Goldman bracket, and the Turaev cobracket and on 
cobordisms of free knots. What is a {\it free knot}? In this book's terminology it is a flat virtual knot taken up to virtualization consisting in allowing one to exchange adjacent virtual and real crossings. This frees
up part of the restriction at a vertex. Equivalently, a free knot is represented by a Gauss diagram that has no signs and no arrows, taken up to abstract Reidemeister moves. There is a natural forgetful functor from virtual knots to free knots, and if a virtual knot goes to a non-trivial free knot, then this virtual knot is itself non-trivial and non-classical. The theory of free knots has been pursued with much energy by
Manturov and his collaborators. At first it was thought (a conjecture of Turaev) that free knots were trivial, but Manturov showed, using parity, that this is not the case and that there are non-trivial cobordism classes of free knots.. The analogous virtual knot theory is to study virtual knots (with over and undercrossings) up to change of orientation of the crossings as we desciribe at the beginning of the review, that is to study virtual knots up to virtualization equivalence. The theory of free knots sheds 
light on these questions of virtualization. In studying free knots we are at one of the most important frontiers of this subject and its relationships with quantum topology \cite{SL3}. And finally there is Chapter 9 on Graph Links, a
further combinatorial generalization of the theory that has much promise in illuminating its original
problems.
\smallbreak

We hope that the reader of this review is motivated to delve into the adventure presented by this 
remarkable book.
\bigbreak

\end{document}